\documentclass[11pt, reqno]{amsart}
\usepackage[utf8]{inputenc}
\usepackage{amsmath, amssymb, amsfonts, latexsym, amsthm}
\usepackage{calligra}
\usepackage{tcolorbox}

\newtheorem{thm}{Theorem}[section]

\numberwithin{equation}{section}

\newcommand{\mbb}{\mathbb}
\newcommand{\ra}{\rightarrow}

\newcommand{\pa}{\partial}
\newcommand{\ov}{\overline}

\begin{document}
\title{Quasiconformal variants of the Wong--Rosay theorem}
    \thanks{}
	\author{Steven G. Krantz, Kaushal Verma} 

\address{SGK: Department of Mathematics, Campus Box 1146, Washington University in St. Louis, St. Louis, Missouri 63130, USA}
\email{sk@math.wustl.edu}

\address{KV: Department of Mathematics, Indian Institute of Science, Bengaluru 560012, India}
\email{kverma@iisc.ac.in}

\begin{abstract}
The Wong--Rosay theorem provides a characterization of the unit ball among all strongly pseudoconvex domains in terms of holomorphic automorphism group actions. We explore variants of this theorem in the quasiconformal setting.
\end{abstract}

\maketitle

\section{Introduction}

Let ${\rm Aut}(D)$ denote the group of holomorphic automorphisms of a bounded domain $D \subset \mathbb C^n$. The natural action of 
${\rm Aut}(D)$ on $D$ given by $(f, z) \mapsto f(z)$, where $(f, z) \in {\rm Aut}(D) \times D$, is proper, real analytic, and ${\rm Aut}(D)$ is locally compact in the compact open topology. Further, ${\rm Aut}(D)$ is a finite dimensional real Lie group that is compatible with the compact open topology. A key ingredient that binds the geometric and function-theoretic features in this set-up is Cartan's compactness principle. It states that every sequence $\{\phi_j\} \in {\rm Aut}(D)$ contains a convergent subsequence and the limit map $\phi$ either maps $D$ entirely into $\partial D$, i.e., $\phi(D) \subset \partial D$ or else $\phi \in {\rm Aut}(D)$. For $a \in D$, let $\mathcal O_a = \{f(a): f \in {\rm Aut(D)}\}$ be its orbit. Call $p \in \partial D$, a {\it boundary accumulation point} if $p \in \overline {\mathcal O_a} \cap \partial D$ for some $a \in D$. 

According to the Wong--Rosay theorem (\cite{Ro}, \cite{W}), if $\partial D$ is smooth strongly pseudoconvex near a boundary accumulation point $p$, then $D$ is biholomorphic to the unit ball $\mathbb B^n \subset \mathbb C^n$. Thus, from a purely local assumption about the domain $D$ near $p$, a global conclusion is being drawn. The passage from local to global is typical of all the further developments in this line and two notable results are the Bedford--Pinchuk theorems in \cite{BP1}, \cite{BP2} in which the scaling method plays an important role. The Wong--Rosay paradigm has been further explored by several others and a survey by Isaev--Krantz \cite{IK} provides an overview of this area. The Greene--Krantz conjecture that predicts the geometric nature of the boundary of a domain near boundary accumulation points is a central question that remains open in its complete generality. 

\medskip

The purpose of this note is to explore quasiconformal variants of the Wong--Rosay theorem. The motivation for doing so comes in part from the notion of {\it quasiconformal homogeneity} that goes back to Gehring--Palka \cite{GP1}. 

\medskip

Consider the standard one-point compactification $\widehat {\mathbb R}^n = \mathbb R^n \cup \{ \infty \}$ of $\mathbb R^n$ equipped with the spherical metric. A domain $D \subset \widehat {\mathbb R}^n$ is said to be {\it quasiconformally homogeneous} if, for a given pair of points $x, y \in D$, there is a quasiconformal map $f : D \rightarrow D$ such that $f(x) = y$. While most domains in $\mathbb C^n$ do not have non-trivial holomorphic automorphisms, the situation is quite the opposite in the quasiconformal setting. Indeed, in \cite{GP1}, Gehring--Palka show that {\it every} domain $D \subset \widehat {\mathbb R}^n$ is quasiconformally homogeneous. This can be interpreted to say that every domain $D \subset \widehat {\mathbb R}^n$ admits a non-trivial group of quasiconformal self-maps. The Gehring--Palka construction is rather explicit and it shows that the maximal dilatations of such maps, as the points $x, y$ vary in $D$, are not uniformly bounded.

\medskip

This leads to the related notion of a {\it uniformly quasiconformally homogeneous} domain. A domain $D \subset \widehat {\mathbb R}^n$ is said to be $K$-quasiconformally homogeneous if for all $x, y \in D$, there is a $K$-quasiconformal map $f: D \rightarrow D$ such that $f(x) = y$. If $D$ is
$K$-quasiconformally homogeneous for some $K$, then $D$ is said to be uniformly quasiconformally homogeneous. Not all domains are uniformly quasiconformally homegeneous. Theorem $5.5$ in \cite{GP1} shows that the boundary $\partial D$ of such a domain must satisfy one of the following three conditions: $\partial D$ has at most two points, $\partial D$ is a non-degenerate continuum, or every neighbourhood of each point on $\partial D$ has infinitely many components of $\partial D$. In addition, Theorem $7.4$ of \cite{GP1} shows that the last possibility above cannot occur if $D$ is uniformly quasiconformally homogeneous with respect to a group of such maps. Other variations of the notion of quasiconformally homogeneity can be found in \cite{BCMTW}, \cite{GM}. The survey \cite{BCT} provides an overview of some other aspects of quasiconformal homogeneity.

\medskip

A basic problem that has been studied extensively is that of characterizing the domains $D \subset \widehat {\mathbb R}^n$ that can be mapped quasi-conformally onto the unit ball 
\[
B^n(0, 1) = \{x = (x_1, x_2, \ldots, x_n) \in \mbb R^n: \Vert x \Vert < 1\}
\]
in $\mathbb R^n$. Here, $\Vert x \Vert^2 = x_1^2 + x_2^2 + \ldots + x_n^2$ is the standard euclidean norm on $\mbb R^n$. More details on this can be found in \cite{G1}, \cite{GMP} and \cite{V}. A notion that has been used to study this basic problem, and which is relevant here as well, is that of a {\it quasiconformally collared} boundary point. Before recalling the definition, let
\[
\mathbb H^n = \{x = (x_1, x_2, \ldots, x_n) \in \mathbb R^n: x_n > 0\}
\]
be the upper half-plane in $\mathbb R^n$. For a domain $D \subset \mathbb R^n$ and $p \in \partial D$, recall that $D$ is {\it quasiconformally collared at $p$}, if there is a neigbourhood $U$ of $p$ in $\mathbb R^n$ and a homeomorphism $f$ of $U \cap \overline D$ onto $B^n(0, 1) \cap \overline {\mathbb H}^n$ such that its restriction to $U \cap D$ is quasiconformal. Evidently, $f$ is a quasiconformal surjection of $U \cap D$ onto $B^n(0,1) \cap \mathbb H^n$ and $f(U \cap \partial D) = B^n(0, 1) \cap \partial \mathbb H^n$. By composing with a M\"obius transformation that preserves $B^n(0,1) \cap \mathbb H^n$, we can always assume that $f(p) = 0 \in \partial \mathbb H^n$. 

\medskip

The class of pairs $(D, p)$, with $p \in \partial D$, such that $D$ is quasiconformally collared at $p$ is rather large. As mentioned in \cite{GMP}, it includes domains that look like wedges, cones, or have Lipschitz boundaries near $p$. To clarify, wedges $W(n, \alpha), 0 < \alpha < 2 \pi$, are domains in $\mathbb R^n$ that are defined using cylindrical coordinates in $\mathbb R^2$, i.e.,
\[
W(n, \alpha) = \{(r, \theta, z): r > 0, 0 < \theta < \alpha, z \in \mathbb R^{n-2}\}
\]
and cones $C(n, \alpha), 0 < \alpha < \pi$ (the standard open cone with vertex at the origin and half-vertex angle $\alpha$) are domains in $\mathbb R^n$ defined as
\[
C(n, \alpha) = \{x \in \mathbb R^n: x_n > \vert x \vert \cos \alpha\}.
\]

The first main result of this paper is as follows:

\begin{thm}
    Let $D \subset \mathbb R^n$, $n \ge 2$ be a domain that is quasiconformally collared at $p \in \partial D$. Suppose there is a family of $K$-quasiconformal maps $\phi_j : D \rightarrow D$ and $a \in D$ such that $\phi_j(a) \rightarrow p$. Then $D$ is quasiconformally equivalent to the unit ball $B^n(0, 1)$.
\end{thm}

Some remarks are in order. First, the parallels with the Wong--Rosay theorem are clearly seen. From a local hypothesis, namely the existence of a quasiconformally collared boundary point at which an orbit accumulates, a global conclusion, namely the quasiconformal equivalence of $D$ with the ball, is drawn. 

\medskip

Second, the hypothesis that a given point $a \in D$ is mapped by quasiconformal maps to points that converge to $p$ is evidently much weaker than $D$ being quasiconformally homogeneous. 

\medskip

Third, Gehring (\cite{G2}) has constructed bounded Jordan domains $D \subset \mathbb R^n$, $n \ge 3$, such that $D$ has a tangent plane $T$ at all points of $\partial D$ such that $T$ is continuous at all points except one, and such that $D$ is {\it not} quasiconformally equivalent to $B^n(0,1)$. Note that $\partial D$ is $C^1$-smooth at points where $T$ is continuous and hence it is quasiconformally collared there. By \cite{GP1}, $D$ is also quasiconformally homogeneous and therefore, quasiconformal orbits of points in $D$ must accumulate everywhere on $\partial D$. Nevertheless, Theorem $1.2$ shows that none of the orbits that accumulate at points where $T$ is continuous can have uniformly bounded dilatations. If they did, $D$ would be quasiconformally equivalent to $B^n(0,1)$, which it is not by construction.

\medskip

A second variant of the Wong--Rosay theorem can be described in the context of quasiconformal harmonic maps. To describe it, let $D \subset \mathbb C$ be a domain. A complex-valued map $f: D \rightarrow \mathbb C$ is said to be {\it harmonic} if $f = u + iv$ where both $u, v$ are harmonic on $D$. The real and imaginary parts $u, v$ are not required to satisfy the Cauchy--Riemann conditions. A clear exposition of the theory of such maps can be found in Duren's book \cite{Du}. If $D$ is simply connected, it is known that there are holomorphic functions $g, h$ on $D$ such that $f = g + \overline h$. While such a representation always holds locally, the presence of non-trivial periods is an obstruction for its global validity in multiply connected domains.

\medskip

For an orientation-preserving diffeomorphism $f: D \rightarrow f(D) \subset \mbb C$ to be quasiconformal, the {\it complex dilatation} $\mu_f = f_{\overline z}/f_z$ must satisfy $0 \le \vert \mu_f \vert \le k < 1$ in $D$. For harmonic maps, the {\it second complex dilatation} $\omega_f = \overline {f_{\overline z}}/f_z$ is more relevant. Note that $\vert \mu_f \vert = \vert \omega_f \vert$. For a quasiconformal harmonic diffeomorphism $f$, $\vert \omega_f \vert = \vert \mu_f \vert \le k < 1$. By writing $f = g + \overline h$ locally, it is seen that $\nu_f = h_z/g_z$ is a bounded holomorphic function on $D$. In addition, the Jacobian $J_f$, which can always be expressed as $J_f = \vert f_z \vert^2 - \vert f_{\overline z} \vert^2$, is seen to satisfy the inequality $J_f \ge \vert g_z \vert^2(1 - k^2)$. 

\medskip

The Rad\'o--Kneser--Choquet theorem (see \cite{Du}) asserts that the harmonic extension of an orientation-preserving homeomorphism of the unit circle 
continues to a harmonic diffeomorphism of the open unit disc that extends continuously to its boundary. Generalizations of this theorem are known for bounded convex domains with Jordan boundaries. This gives a natural way to create harmonic self-maps of a given convex domain with Jordan boundaries. There are versions of this theorem for finitely connected domains also each of whose boundary components is a Jordan curve--see \cite{DH} for example. This motivates the next statement.

\medskip

Our second main result is this:

\begin{thm}
    Let $D \subset \mathbb C$ be a bounded domain and $p \in \partial D$ a $C^2$-smooth boundary point. Suppose that there exist a point $a \in D$ and a family of $k$-quasiconformal harmonic diffeomorphisms $\phi_j : D \rightarrow D$ such that $\phi_j(a)$ converges to $p$. Then there exists a $k$-quasiconformal harmonic diffeomorphism $\phi: D \rightarrow \mathbb H^2$.
\end{thm}

\section{Proof of Theorem 1.1}

Some basic compactness results for families of quasiconformal mappings are required, and we rely on the expositions in \cite{GMP}, \cite{L} and \cite{V}. 

\medskip

\noindent {\it Result 1}: Let $\mathcal F$ be a family of $K$-quasiconformal mappings defined on a domain $D \subset \mathbb R^n \cup \{\infty\}$. If every $f \in \mathcal F$ omits two distinct fixed values, then $\mathcal F$ is equicontinuous in $D$.

\medskip

Viewing each member $f \in \mathcal F$ as a map with values in $\mathbb R^n \cup \{\infty\}$ endowed with the spherical metric, the Arzela--Ascoli theorem shows that $\mathcal F$ is normal. When the family $\mathcal F$ takes values in a varying set of domains, a more precise convergence theorem holds. Associated with a family of domains $D_j \subset \mathbb R^n \cup \{\infty\}$ is its {\it kernel} which will be denoted by 
${\rm ker}(D_j)$ and defined as follows:
\[
{\rm ker}(D_j) = \bigcup_{\nu \ge 1} {\rm int} \left( \bigcap_{\mu \ge \nu} D_{\mu} \right).
\]
In other words, a point $x$ belongs to ${\rm ker}(D_j)$ precisely when it has an open neighbourhood that belongs to all the $D_j$'s for $j$ large. Therefore, a compact subset of ${\rm ker}(D_j)$ is compactly contained in all the $D_j$'s for $j$ large. Note that ${\rm ker}(D_j)$ is open by definition but may not always be connected even though each $D_j$ is a domain. The complement of a set $A \in \mathbb R^n \cup \{\infty\}$ will be denoted by $A^c$.

\medskip

\noindent {\it Result 2}: Let $\{D_j\}$ be a sequence of domains in $\mathbb R^n \cup \{\infty\}$, $n \ge 2$, whose kernel contains at least two points in its complement. Let $D \subset {\rm ker}(D_j)$ be a component. Suppose that $f_j : D_j \rightarrow D'_j$ is a family of $K$-quasiconformal mappings that converges locally uniformly to $f$ in $D$. Then one of the two possibilities must hold. Either $f$ is a constant mapping of $D$ that takes its value in $({\rm ker}(D'_j) \cup {\rm ker}((D'_j)^c))^c$ or $f : D \rightarrow D'$ is a $K$-quasiconformal mapping onto a component $D'$ of ${\rm ker}(D'_j)$. In the latter case, the inverses $f_j^{-1}$ converge locally uniformly to $f^{-1}$ in $D'$. 
\medskip \\

\noindent {\it Proof of Theorem 1.1}: The domain $D$ is given to be quasiconformally collared at $p \in \pa D$ via a quasiconformal homeomorphism $f$ that is defined in a one-sided neighbourhood $U \cap D$ of $p$, extends continuously up to $\pa D$, and straightens it near $p$. Thus it is possible to pick a $q \in \pa D$, $q \not= p$ and $q$ near $p$. The family of $K$-quasiconformal maps $\phi_j: D \ra D$ is normal since each $\phi_j$ misses the pair of points $p, q$. The limit $\phi_{\infty}$ of each locally uniformly convergent subsequence of $\{\phi_j\}$ cannot be a homeomorphism since $\phi_j(a) \rightarrow \phi_{\infty}(a) = p \in \partial D$. Hence, $\phi_{\infty}$ must be the constant map and $\phi_{\infty}(z) \equiv p$ for $z \in D$. It follows that the entire sequence $\{\phi_j\}$ locally uniformly converges on $D$ to the constant limit $\phi_{\infty} \equiv p$. 

\medskip

An observation that is also made in the proof of the Wong--Rosay theorem at this stage is the fact that $\phi_{\infty} \equiv p$ on $D$ is sufficient to conclude that $D$ is simply connected. Indeed, a loop $\gamma \subset D$ is contractible precisely when $\phi_j(\gamma)$ is so. But then  $f \circ \phi_j(\gamma) \subset B^n(0,1) \cap \mathbb H^n$ for all large $j$ which implies that $f \circ \phi_j(\gamma)$ and hence $\phi_j(\gamma)$ is contractible for large $j$. 

\medskip

Moving on with the proof, to get a non-degenerate limit mapping, it is necessary to modify the $\phi_j$s appropriately.

\medskip

Note that $f \circ \phi_j(a) \rightarrow p_0 \in 0 \in \partial \mathbb H^n$. Also, for each compact $C \subset D$, the maps $f \circ \phi_j$ are well defined on $C$ for large $j$ and map it into $B^n(0, 1) \cap \mathbb H^n$. Choose a family of M\"obius maps $\psi_j \in {\rm Aut}(\mathbb H^n)$ such that $\psi_j \circ f \circ \phi(a) = b = (0, 0, \ldots, 1)\in \mathbb H^n$ for all $j$. Since the $\psi_j^{-1}$'s converge locally uniformly on $\mathbb H^n$ to $0 \in \partial \mathbb H^n$, the domains $\psi_j^{-1}(B^n(0,1) \cap \mathbb H^n)$ fill up all of $\mathbb H^n$. 

\medskip

Let $f$ be $K'$-quasiconformal. With $D_1 \subset D_2 \subset \ldots \subset D$ an exhaustion of $D$ by an increasing sequence of compactly contained domains, let $D'_{lj} = \psi_j \circ f \circ \phi(D_l) \subset \mathbb H^n$. Then, for each $l$ there is a $j_l$ depending on $l$ such that 
\[
T_j = \psi_j \circ f \circ \phi: D_l \ra D'_{lj} \subset \mathbb H^n
\]
is well defined for $j \ge j_l$ and they form a family of quasiconformal maps whose maximal dilatation is at most $KK'$. Note that $T_j(a) = b$ for all $j$ and ${\rm ker}(D_l) = D$. The family $\mathcal T = \{T_j\}$ is normal on $D$ since each member takes values in $\mathbb H^n$ and therefore misses at least two fixed distinct points. 

\medskip

Let $H_1 \subset H_2 \subset \ldots \subset \mathbb H^n$ be an exhaustion of $\mathbb H^n$ by an increasing sequence of compactly contained domains. Let $H'_{lj} = \phi^{-1} \circ f^{-1} \circ \psi_j^{-1}(H_l) \subset D$. Then, for each $l$ there is a $j_l$ depending on $l$ such that
\[
S_j = \phi^{-1} \circ f^{-1} \circ \psi_j^{-1}: H_l \ra H'_{lj} \subset D
\]
is well defined for $j \ge j_l$ and they form a family of quasiconformal maps whose maximal dilatation is at most $KK'$. Note that $S_j(b) = a$ for all $j$ and ${\rm ker}(H_l) = \mathbb H^n$. The family of inverses $\mathcal S = \{S_j\}$ is normal on $\mathbb H^n$ since each member takes values in $D$ and therefore misses at least two fixed distinct points.  

\medskip

Fix $H_{l_0} \subset \mathbb H^n$ and let $\widetilde S$ be a uniform limit on $H_{l_0}$ of a subsequence $\{S_{j_i}\}$. Then, either $\widetilde S$ is a constant map or a homeomorphism onto its image which lies in $\overline D$. The former case cannot hold, for if so, then $\widetilde S \equiv a$ on $H_{l_0}$. Then, for all large $j_i$, the images $S_{j_i}(H_{l_0})$ are contained in a small ball around $a$. The maps $T_{j_i}$ are defined on this ball for all large $j_i$, and by normality, admit a further subsequence that converges uniformly on it. Let $\widetilde T$ be the limit of $T_{j_i}$. Retaining the same notation for such a subsequence, note that $T_{j_i} \circ S_{j_i} = {\rm id}$ on $H_{l_0}$. Taking limits, $\widetilde T \circ \widetilde S = {\rm id}$ and this shows that $\widetilde S$ cannot be a constant map. Therefore, $\widetilde S: H_{l_0} \ra \widetilde S(H_{l_0}) \subset \overline D$ is a homeomorphism. Now, the range $\widetilde S(H_{l_0})$ equals ${\rm ker}(S_{j_i}(H_{l_0})$. If $s \in \widetilde S(H_{l_0}) \cap \pa D$, it would mean that there is an open neighbourhood of $s$ that is contained in $S_{j_i}(H_{l_0})$ for large $j_i$. This is not possible since the $S_{j_i}$s take values in $D$. Thus, every normal limit of $\mathcal S$ is a homeomorphism onto its image that lies in $D$. The same holds for the family $\mathcal T$. 

\medskip

Now, by passing to the diagonal subsequence in $\mathcal S$, we get a sequence that converges on all the $C_j$s. For brevity, retain the same notation $S_j$ for this subsequence and let $S_{\infty}$ be its uniform limit on all compact sets in $\mathbb H^n$. Then, as explained above, $S_{\infty}$ is a homeomorphism onto its image that lies in $D$. It remains to show that $S_{\infty}$ is surjective. Suppose that $t \in D \cap \pa S_{\infty}(\mathbb H^n)$. Let $V \subset D$ be a small neighbourhood of $t$. Corresponding to the diagonal subsequence $S_j$ above, consider the inverses $T_j$ and pass to a further subsequence, if required, so that 
$\{T_{i_j}\}$ converges uniformly on $V$. As explained above, the limit of $T_{i_j}$, say $T_V$ is a homeomorphism onto its image that lies in $\mathbb H^n$. Since $S_{i_j} \circ T_{i_j} = {\rm id}$ on $V$, the local uniform convergence of the $S_{i_j}$'s allows us to pass to limits to get $S_{\infty} \circ T_V = {\rm id}$ on $V$. This shows that $S_{\infty}$ is surjective.

\section{Proof of Theorem 1.2}

\noindent The starting point of the proof is the same as before. Pick $q \in \partial D$, $q \not= p$ and $q$ near $p$ and observe that the family of $k$-quasiconformal harmonic diffeomorphisms $\phi_j : D \rightarrow D$ is normal since each $\phi_j$ misses the pair of points $p, q$. The limit $\phi_{\infty}$ of each locally uniformly convergent subsequence of $\{\phi_j\}$ cannot be a homeomorphism since $\phi_j(a) \rightarrow \phi_{\infty}(a) = p \in \partial D$. Hence, $\phi_{\infty}$ must be the constant map and $\phi_{\infty}(z) \equiv p$ for $z \in D$. It follows that the entire sequence $\{\phi_j\}$ converges locally uniformly on $D$ to the constant limit $\phi_{\infty} \equiv p$. As before, this implies that $D$ is simply connected. 

\medskip

Let $\rho$ be a $C^2$-smooth defining function for $\partial D$ near $p$ such that, in a neighbourhood $U$ of $p$, the sub-level set $U \cap D = \{\rho < 0\}$ and $d\rho \not= 0$ on $U \cap \partial D$ . Apply the scaling method now with respect to the orbit sequence $\phi_j(a)$. For brevity, write $\phi_j(a) = p_j$. The affine maps
\[
T_j(z) = \frac{z - p_j}{-\rho(p_j)}
\]
satisfy $T_j(p_j) = 0$ and every compact subset in the plane is eventually contained in $T_j(U)$ for $j$ large since $\rho(p_j) \rightarrow 0$. Fix an arbitrarily large ball $B(0, R)$. The normalized functions 
\begin{align*}
\frac{1}{-\rho(p_j)} \rho \circ T_j^{-1}(z) &= \frac{1}{-\rho(p_j)} \rho(p_j - \rho(p_j))\\
                                            &= \frac{1}{-\rho(p_j)} \left(\rho(p_j) - 2 \rho(p_j) \Re(\partial \rho(p_j) z) + (\rho(p_j))^2 O(1)  \right)\\
                                            &= -1 + 2 \Re(\partial \rho(p_j) z) - \rho(p_j) O(1)
\end{align*}
are then well defined on $B(0, R)$ for large $j$ and serve as a $C^2$-smooth defining function for the domain $T_j(U \cap D)$. The $O(1)$ term is locally uniformly bounded on $\mathbb C$. Let
\[
\mathbb H^{\mp} = \{z \in \mathbb C: {\pm}(-1 + 2 \Re(\partial \rho(p) z)) < 0\}.
\]
It can be seen that the domains $T_j(D)$ and $T_j(U \cap D)$ both converge in the Hausdorff sense on all compact subsets of the plane to the limiting half-plane $\mathbb H^-$. It follows that ${\rm ker}(T_j(D)) = \mathbb H^-$ and ${\rm ker}((T_j(D))^c) = \mathbb H^+$.

\medskip

Thus, we are led to study the family $\mathcal F$ consisting of $F_j = T_j \circ \phi_j: D \rightarrow T_j(D)$. While the class of harmonic maps is preserved under pre-composition by holomorphic maps, this is not true in general for post-compositions. The redeeming feature here is that $T_j$ is an affine holomorphic map. Since the factor $\rho(p_j)$ is real, the chain rule shows that
\[
\frac{\overline {(F_j)_{\ov z}}}{(F_j)_z} = \frac{\overline {(\phi_j)_{\ov z}}}{(\phi_j)_z}
\]
and hence each $F_j$ is a quasiconformal harmonic diffeomorphism onto its image with uniformly bounded second complex dilatation.

\medskip

To show that $\mathcal F$ is normal in $D$, fix a compact $G \subset D$. First, note that the maximal dilatation of each $F_j$ is at most $K = (1+k)/(1-k)$. Second, as observed above, $\phi_j(G) \subset U \cap D$ for all large $j$. Therefore $F_j = T_j \circ \phi_j$, when restricted to $G$, is a family of quasiconformal maps, with uniformly bounded maximal dilatations, that eventually takes values in $T_j(U \cap D)$. This is a domain that is described in part by $(-\rho(p_j))^{-1} \rho \circ T_j^{-1}$ which is approximately a line as seen above. It follows that each $F_j$, when restricted to $G$, misses a suitably chosen open neighbourhood in $\mathbb H^+$. It is therefore equicontinuous and hence normal when viewed as maps into $\mathbb R^2 \cup \{\infty\}$. To conclude, $F_j = T_j \circ \phi_j$ converges locally uniformly in $D$ after passing to a suitable subsequence. The limit cannot be a constant map, for in this case, the limiting constant value must lie on the line
\[
-1 + 2 \Re(\partial \rho(p)z) = 0
\]
and this cannot happen since $T_j \circ \phi_j(a) = 0$. The limit map $F_{\infty}$ is then a $k$-quasiconformal map from $D$ onto $\mathbb H^{-}$. It is also harmonic being the uniform limit of such maps. An affine holomorphic map transforms $\mathbb H^-$ onto $\mathbb H^2$. The desired $k$-quasiconformal diffeomorphism from $D$ onto $\mathbb H^2$ is obtained by post-composing $F_{\infty}$ with this affine holomorphic map.


\section{Concluding Remarks}

\noindent The theorem of Bun Wong and Rosay has been quite influential in
complex function theory for the past $45$ years.  It has led to many
new insights about automorphism group actions.  We mention particularly
the Greene-Krantz conjecture.  

\medskip

It is appropriate now to generalize the idea of the Wong--Rosay theorem
to other contexts.  The idea of quasi-conformal mapping is certainly
one natural candidate for such a generalization.   We hope in future
work to explore an axiomatic approach to this set of ideas.




\begin{thebibliography}{BCMTW}

\bibitem{BP1}
Bedford, E.; Pinchuk, S. I.:
\textit{Convex domains with noncompact groups of automorphisms}. (Russian) Mat. Sb. 185 (1994), no. 5, 3–26; translation in Russian Acad. Sci. Sb. Math. 82 (1995), no. 1, 1–20.

\bibitem{BP2}
 Bedford, Eric; Pinchuk, Sergey:
 \textit{Domains in $\mathbb C^2$ with noncompact automorphism groups}. Indiana Univ. Math. J. 47 (1998), no. 1, 199–222.

\bibitem{BCT}
Bonfert-Taylor, Petra; Canary, Richard D.; Martin, Gaven; Taylor, Edward C.:
\textit{Quasiconformal homogeneity after Gehring and Palka}, 
Comput. Methods Funct. Theory {\bf 14}  (2014), no. 2-3, 417–430.

\bibitem{BCMTW}
Bonfert-Taylor, Petra; Canary, Richard D.; Martin, Gaven; Taylor, Edward C. ;Wolf, Michael:
\textit{Ambient quasiconformal homogeneity of planar domains}, (English summary) Ann. Acad. Sci. Fenn. Math. {\bf 35} (2010), no. 1, 275–283.

\bibitem{Du}
Duren, Peter:
\textit{Harmonic mappings in the plane}, {\bf 156}, Cambridge University Press. 

\bibitem{DH}
Duren, P.; Hengartner, W.: 
\textit{Harmonic mappings of multiply connected domains}, Pacific J. Math. {\bf 180} (1997), 201--220.


\bibitem{GKK}
Gaussier, Herv\'e; Kim, Kang-Tae; Krantz, Steven G.:
\textit{A note on the Wong-Rosay theorem in complex manifolds}, Complex Var. Theory Appl. {\bf 47} (2002), no. 9, 761–768.

\bibitem{G1}
Gehring, F. W.:
\textit{Extension theorems for quasiconformal mappings in $n$-space}, J. Anal. Math. {\bf 19} (1967), 149-169.

\bibitem{G2}
Gehring, F. W.:
\textit{A remark on domains quasiconformally equivalent to a ball}, Ann. Acad. Sci. Fenn. Ser. A I Math. 2 (1976), 147–155

\bibitem{GMP}
Gehring, F. W.; Martin, G.; Palka, B.: 
\textit{An introduction to the theory of higher dimensional quasiconformal mappings}, Math surveys and monographs {\bf 216}, Amer. Math Soc.

\bibitem{GP1}
Gehring, F. W.; Palka, B. P.:
\textit{Quasiconformally Homogeneous Domains}, J. Anal. Math. {\bf 30} (1976), 172–199.

\bibitem{GM}
Gong, Jianhua; Martin, Gaven:
\textit{Aspects of quasiconformal homogeneity} (English summary)
New Zealand J. Math. {\bf 39} (2009), 117–132.

\bibitem{IK}
Isaev, A. V.; Krantz, S. G.:
\textit{Domains with non-compact automorphism group: a survey}. Adv. Math. 146 (1999), no. 1, 1–38. 

\bibitem{KK}
Kim, Kang-Tae; Krantz, Steven G.:
\textit{A Kobayashi metric version of Bun Wong's theorem}, Complex Var. Elliptic Equ. {\bf 54} (2009), no. 3-4, 355–369.

\bibitem{L}
Lehto, O.:
\textit{Univalent functions and Teichm\"uller theory}, Graduate Texts in Mathematics, 109, Springer.

\bibitem{Ro}
Rosay, Jean-Pierre:
\textit{Sur une caract\'erisation de la boule parmi les domaines de $\mathbb C^n$ par son groupe d'automorphismes}. (French) Ann. Inst. Fourier (Grenoble) {\bf 29} (1979), no. 4, ix, 91–97.

 \bibitem{V}
 V\"ais\"al\"a, J.:
 \textit{Lectures on $n$-dimensional quasiconformal mappings}, Lecture Notes in Mathematics, 229, Springer.

 \bibitem{W}
 Wong, B.:
 \textit{Characterization of the unit ball in $\mathbb C^n$ by its automorphism group}. Invent. Math. {\bf 41} (1977), no. 3, 253–257.
 
\end{thebibliography}
\end{document}